\documentclass[journal,comsoc]{IEEEtran}

\usepackage[T1]{fontenc}

\usepackage{amsmath}
\usepackage[cmintegrals]{newtxmath}
\usepackage{algorithmic}
\usepackage{algorithm}
\usepackage{array}
\usepackage[caption=false,font=normalsize,labelfont=sf,textfont=sf]{subfig}
\usepackage{textcomp}
\usepackage{stfloats}
\usepackage{url}
\usepackage{verbatim}
\usepackage{graphicx}
\usepackage{cite}
\usepackage{xcolor}
\usepackage{soul}
\usepackage{url}
\usepackage{hyperref}

\begin{document}

\title{An Analytical Model for Performance Estimation in High-Capacity IMDD Systems}

\author{Giuseppe~Rizzelli, Pablo~Torres-Ferrera, Fabrizio~Forghieri,~\IEEEmembership{Fellow,~IEEE} and~Roberto~Gaudino,~\IEEEmembership{Senior Member,~IEEE}
\thanks{G. Rizzelli is with the LINKS Foundation, Torino, Italy, e-mail: giuseppe.rizzelli@linksfoundation.com.}
\thanks{P. Torres-Ferrera is with the Department of Engineering, University of Cambridge, Cambridge, United Kingdom. Prior, he was with the Department of Electronics and Telecommunications, Politecnico di Torino, Torino, Italy.}
\thanks{F. Forghieri is with CISCO Photonics Italy, Vimercate, Italy.}
\thanks{R. Gaudino is with the Department
of Electronics and Telecommunications, Politecnico di Torino, Torino,
Italy.}}
\date{March 2023}

\maketitle

\begin{abstract}
In this paper, we propose an analytical model to estimate the signal-to-noise ratio (SNR) at the output of an adaptive equalizer in intensity modulation and direct detection (IMDD) optical transmission systems affected by shot noise, thermal noise, relative intensity noise (RIN), chromatic dispersion (CD) and bandwidth limitations. We develop the model as an extension of a previously presented one, and then we test its accuracy by sweeping the main parameters of a 4-PAM-based communication system such as RIN coefficient, extinction ratio, CD coefficient and equalizer memory. Our findings show a remarkable agreement between time-domain simulations and analytical results, with SNR discrepancies below 0.1 dB in most cases, for both feed-forward and decision-feedback equalization. We consider that the proposed model is a powerful tool for the numerical design of strongly band-limited IMDD systems using receiver equalization, as it happens in most of modern and future M-PAM solutions for short reach and access systems.
\end{abstract}

\begin{IEEEkeywords}
Instensity Modulation, Direct Detection, Optical Communications, Performance Modelling.
\end{IEEEkeywords}

\IEEEpeerreviewmaketitle

\section{Introduction}
\label{introduction}
\IEEEPARstart{T}{he} capacity requirements of short-reach and access networks is surpassing those of long-reach networks. Although coherent detection (CoD) is making its way into shorter transmission reach segments \cite{400G_DCI, scaling_laws, OFC2022, CPON_Faruk, overview, JLT_MMF}, IMDD systems are still the solution of choice in most of the metro edge and access networks applications, as well as in many intra- and inter- data center interconnects (DCI) \cite{overview, IMDD_PS, IMDD_200G, JLT_pon, IMDD_MLSE, JLT2021_MMF_SWDM}. Regardless of the specific transmission scenario, IMDD-based options are increasingly relying on more advanced modulation formats than the traditional on-off keying (OOK) in order to improve legacy systems performance and thus extend their lifespan in an increasingly data-hungry society. For instance, the 400GBASE-LR4 and the 400GBASE-SR standards have defined solutions for SMF-based short-reach and VCSEL+MMF-based DCI links, respectively, relying on 4-PAM modulation to reach up to 100 $Gbps/\lambda$.

As it is the case for the design or analysis of any transmission system, the possibility to numerically predict the performance of an IMDD system can be of great help in the effort to efficiently plan and optimize network capacity at the physical layer. This is usually achieved through computationally intensive simulations in the time domain that allow to estimate, for instance, the bit error ratio (BER) of the system for a given set of fixed parameters. However, such a time consuming approach can be highly impractical when the optimization of a great number of parameters is required as, for instance, in a multi-dimensional Monte Carlo analysis \cite{Cantono} or in a statistical analysis \cite{JLT2021_MMF_SWDM}. In a recent paper \cite{Coh_Fischer} we have presented an analytical tool for performance prediction applied to a coherent system considering dual polarization transmission, a linear channel described by a [2x2] frequency dependent transfer function, colored additive Gaussian noise at the receiver and a full digital signal processing (DSP) implementation at the receiver. The proposed model allows for a 333 times reduction in computation time with respect to a full time domain simulation generating a single output in less than 0.05 seconds on a standard commercial laptop. In this new contribution, we propose a similar approach to compute the electrical SNR at the equalizer output in an IMDD scheme with M-PAM modulation, considering different sources of impairments such as shot noise due to photodection, thermal noise of the transimpedence amplifiers (TIA) at the receiver, RIN associated to the transmitter laser, fiber chromatic dispersion and equalizer implementation. Thus, compared to existing tools, our approach enables fast and accurate performance prediction of modern IMDD systems affected by a broad variety of impairments. Differently from the coherent case, where propagation can be described through an additive white Gaussian noise (AWGN) channel on the optical field, in the IMDD case not all the noise contributions are signal independent. Specifically, the variance of the shot noise is proportional to the instantaneous power of the modulated signal, whereas the RIN scales with the square of this instantaneous optical power. As a time dependent power evolution cannot be taken into account analytically, in our model we define the noise variance proportional to the average optical power (squared for RIN). This approximation is made on the grounds that the equalizer filter operates over a number of received symbols (the equalizer memory), effectively correlating among them and producing a sort of averaging effect of the instantaneous noise (shot or RIN) contributions in time. Thus, we do not expect any estimation error originating from this assumption, at least for band-limited conditions where the effect of the equalizer is more pronounced.

As our coherent analytical model \cite{Coh_Fischer}, the IMDD model proposed here is based on the work by Robert Fischer \cite{Fischer}. The novelty of our approach is in the extension of the tool presented in \cite{Fischer} to include, not only the effect of power independent AWGN such as the TIA thermal noise, but also other noise sources instantaneously proportional to the useful signal power. Moreover, we further refine our analytical model to also take into account the effect of CD and of additional shot noise due to an avalanche photodiode (APD). After introducing all the features of our tool, we compare it against full time domain simulations based on error counting and show the estimation error in terms of both SNR and BER at the equalizer output.

The remainder of this manuscript is then organized as follows: in Section \ref{model} we describe the analytical model, defining the assumptions and showing the main equations. In Section \ref{validity} we show the main results of the comparison of the analytical model performance with a comprehensive set of time domain simulations, highlighting its accuracy in a 4-PAM-based transmission system under several different conditions. In Section \ref{pon} we focus on a specific application scenario and apply our model to a  passive optical networks (PON) in O-band, including CD and APD-based detection. Lastly, in Section \ref{conclusion} we summarize the main results and draw some conclusion.

\section{The proposed analytical model}
\label{model}

\begin{figure}[h]
\centering
\includegraphics[width=1\linewidth]{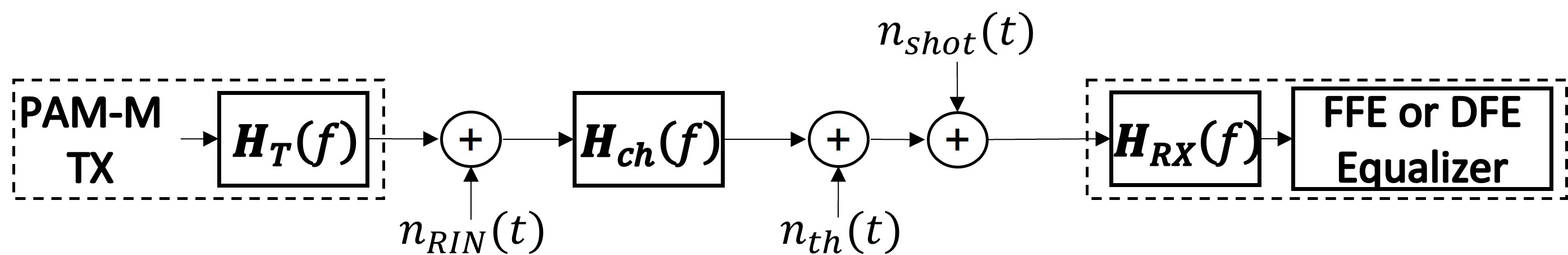}
\caption{Block diagram of the system under investigation with noise sources.}
\label{fig:setup}
\end{figure}

A simplified block diagram of the system under investigation is depicted in Fig. \ref{fig:setup}. It is composed of a transmitter that generates the M-PAM modulated signal and also includes the shaping filtering stage indicated by the transfer function $H_T(f)$. For simplicity, we will show numerical examples assuming a rectangular shape for the transmitted symbols in the time domain, but any other shape can also be handled analytically. Moreover, although our approach can be applied to any M-PAM format, we will mainly focus on 4-PAM modulation. We assume that the RIN at the transmitter can be modeled as a noise source with variance $\sigma^2_{RIN}=k_{RIN}\cdot{P_{TX}^2(t)}$ proportional to the average transmitted power squared, added at the transmitter before the signal is filtered by a linear channel with end-to-end transfer function $H_{ch}(f)$. The channel can be described by a generic frequency response of any type, but in this work we use a supergaussian profile of variable order. Please note that $H_{ch}(f)$ also acts on the power spectral density (PSD) of the RIN. At the receiver side, the filtered signal receives the contributions of thermal and shot noise with PSD $\sigma^2_{th}=N_0$ and $\sigma^2_{shot}=k_{shot}\cdot{P_{RX}(t)}$, respectively. The receiver is equipped with two possible Minimum Mean-Square Error (MMSE)-based equalization schemes, a feedforward equalizer (FFE) and a decision feedback equalizer (DFE). Before the equalizer, a receiver filter with generic transfer function $H_{RX}(f)$ can also be taken into account. In the time domain simulator the BER is then computed through an error counting technique and the SNR calculated as the ratio of the average energy of the signal to the mean square error at the output of the equalizer.

In \cite{Fischer} an analytical model was developed for the estimation of the SNR at the equalizer output, given a generic linear transfer function $H_{ch}(f)$, under two main assumptions: i) the channel is AWGN and ii) the equalizer is infinitely long. A more detailed description of the model can be found in \cite{Fischer} or \cite{Coh_Fischer}. Here, we just recall that the spectrally resolved SNR(f) at the equalizer input is computed as:
\begin{equation}\label{eq:SNR_f}
    SNR(f) = \frac{T \cdot P_{TX}\cdot |H_T(f)\cdot H_{ch}(f)|^2}{N_0(f)}
\end{equation}
where $T$ is the symbol period, $P_{TX}$ is the transmitted signal power, $H_{ch}(f)$ is the linear channel transfer function, $H_T(f)$ is the transfer function of the transmitter shaping filter and $N_0(f)$ is the equivalent noise power spectral density at the input of the receiver. Eq. \ref{eq:SNR_f} is taken from \cite{Fischer} and has to be properly modified to be applied to the IMDD case, as explained below. The SNR at the output of the infinitely long equalizer is then computed as:
\begin{equation}\label{eq:SNR_ffe}
    SNR_{FFE} = \frac{1}{T \cdot \int_{-\frac{1}{2T}}^{\frac{1}{2T}} \frac{1}{\overline{SNR}(f)+1} df}-1
\end{equation}
\begin{equation}\label{eq:SNR_dfe}
    SNR_{DFE} = e^{T \cdot \int_{-\frac{1}{2T}}^{\frac{1}{2T}} log[{\overline{SNR}(f)+1]} df}-1
\end{equation}
respectively for FFE and DFE equalization, where $\overline{SNR}(f)$ is the spectral $SNR(f)$ folded on a bandwidth equal to the symbol rate, defined as:
\begin{equation} \label{eq:SNR_folded}
    \overline{SNR}(f) =  \sum_{\mu} SNR(f - \frac{\mu}{T})
\end{equation}
where the integer $\mu$ indicates the number of foldings.

From Eq. \ref{eq:SNR_f} it is clear that the original model in \cite{Fischer} assumes only signal-independent AWGN at the receiver. In our system, on the other hand, we have two noise sources that depend on the instantaneous power of the signal. As we cannot account for the power dependence analytically, we introduce here the first approximation of our approach: the RIN and shot noise variance are proportional to the average optical power of the useful signal and their PSD can thus be written as $PSD_{RIN}=k_{RIN}\cdot\overline{P_{TX}^2}$ and $PSD_{shot}=k_{shot}\cdot\overline{P_{RX}}$, where $k_{RIN}=RIN_{coeff}/2$ and $k_{shot}=2G^2FqR$ are the proportionality factors for RIN and shot noise, respectively, $RIN_{coeff}$ is the RIN coefficient expressed in $1/Hz$, $G$ is the photodetector gain (when considering avalanche photodetection), $F$ is the photodetector excess noise figure, $R$ is the photodiode responsivity and $q$ is the electron charge. Moreover, since modulation is applied to the optical power, particular care has to be taken in the conversion from optical to electrical dB. For a 4-PAM modulated signal we have:
\begin{equation}\label{eq:Ptx}
    P_{TX}(t) = \overline{P_{TX}} + \frac{OMA_{TX}^{outer}}{6} \cdot \sum_{k=- \infty}^{+ \infty} \alpha_k \cdot s(t-kT)
\end{equation}
where $OMA_{TX}^{outer}$ is the outer optical modulation amplitude (OMA), $\alpha_k$ is a random variable taking one of the 4-PAM levels: $[-3,-1,1,3]$  and $s(t-kT)$ is the pulse shape. The PSD of the transmitted signal is thus (neglecting the irrelevant DC component in $f=0$):
\begin{equation}\label{eq:Ptx}
    PSD_{TX}(f) = \frac{T\cdot(OMA_{TX}^{outer})^2}{36} \cdot |S(f)|^2 \cdot \sigma_{\alpha_k}^2
\end{equation}
where $\sigma_{\alpha_k}^2=5$ for 4-PAM and $|S(f)|^2=|H_T(f)|^2$ is the spectral shape of the signal. Thus, for a 4-PAM modulated signal we can write:
\begin{equation}\label{eq:SNR_f_2}
    SNR(f) = \frac{5}{36} \cdot \frac{T \cdot (OMA_{TX}^{outer})^2\cdot |H_T(f)|^2 \cdot |H_{ch}(f)|^2}{PSD_{noise}(f)}
\end{equation}
where $PSD_{noise}(f)=PSD_{th}(f)+PSD_{shot}(f)+PSD_{RIN}(f)\cdot|H_{ch}(f)|^2$ is the sum of the PSD of each noise contribution, $PSD_{th}=N_0/2$, and we have highlighted that each PSD can be colored and that the RIN PSD is also filtered by $H_{ch}(f)$. In the simulator, we define $\sigma^2_{th}=N0 \cdot f_s/2$, $\sigma^2_{shot}=2G^2FqRP_{RX}(t) \cdot f_s/2$ and $\sigma^2_{RIN}=RIN_{coeff}P_{TX}^2(t) \cdot f_s/2$, where $f_s$ is the simulator sampling frequency

\begin{figure}[h]
\centering
\includegraphics[width=1\linewidth]{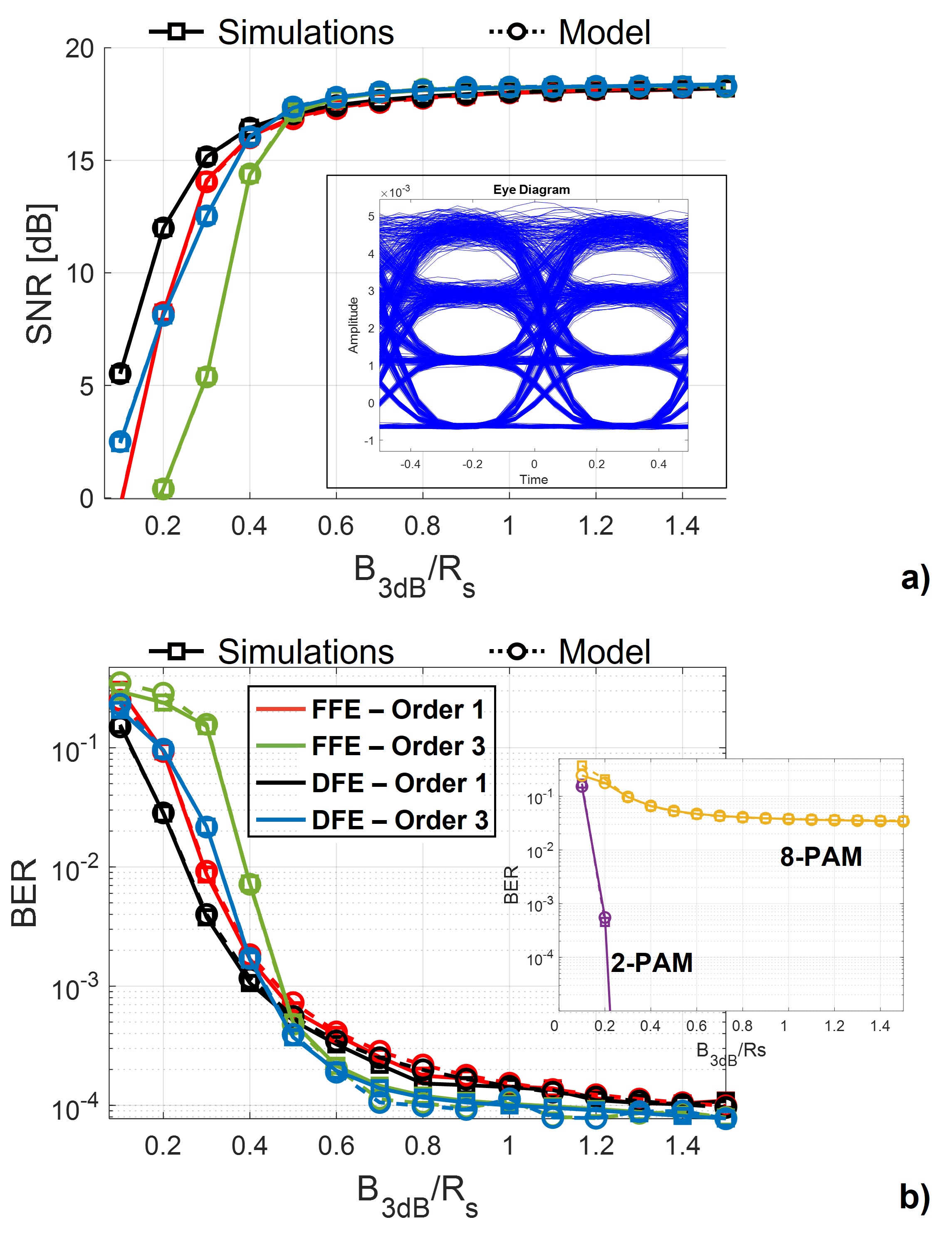}
\caption{a) SNR and b) BER obtained through time domain simulations (solid, squares) and through the proposed analytical model (dashed, circles) as a function of the ratio between the supergaussian filter 3 dB bandwidth $B_{3dB}$ and the symbol rate $R_s$ for a 25 GBaud 4-PAM system, using both FFE (red and green) and DFE (black and blue) equalization. The supergaussian filter order is 1 for black and red curves and 3 for blue and green curves. The inset in a) shows the eye diagram of the 4-PAM signal after the channel filtering when $B_{3dB}/R_s = 0.8$. The inset in b) shows results for FFE-based 2-PAM (purple) and 8-PAM (orange) modulation with supergaussian filter order 1. The legend in b) applies to a) as well.}
\label{fig:model_snr}
\end{figure}

As a first example, Fig. \ref{fig:model_snr} shows a comparison of the results obtained with the time domain simulator and with the proposed analytical model in terms of SNR at the equalizer output as a function of the ratio between the supergaussian filter 3 dB bandwidth ($B_{3dB}$) and the symbol rate ($R_s$) for a 25 GBaud 4-PAM system, using both FFE and DFE equalization. The order of the supergaussian filter is set to 1 or 3. The complete set of simulation parameters are shown in Table \ref{table_param}. These are the parameters used throughout the paper unless otherwise specified. For a fairer comparison between the simulator and the model based on the infinitely long equalizer assumption, in this first model performance validation, we set a high number of taps in the simulator equalizers, respectively 200 for the FFE and 30 for the DFE stage (operated in full-training mode). Fig. \ref{fig:model_snr} shows a very good agreement between the two methods with maximum estimation error of the order of 0.05 dB, regardless of the filter order and of the introduced bandwidth limitations, for both equalization schemes.
\begin{table}[!h]
\caption{Simulation parameters.}
\label{table_param}
\centering
\begin{tabular}{c | cc}
\hline
\textbf{Parameter} & Value & Unit\\
\hline
Symbol Rate & 25 & GBaud\\

RIN & -140 & dB/Hz\\

Extintion Ratio & 6 & dB\\

OMA & 0.78 & dBm\\

Responsivity & 1 & A/W\\

$P_{TX}$ & 1 & mW\\

$N_0$ & $2 \cdot 10^{-19}$ & $W^2/Hz$\\

FFE (DFE) Taps & 200 (30) & \\

\hline
\end{tabular}
\end{table}

The inset in Fig. \ref{fig:model_snr}a shows the eye diagram of the 4-PAM signal after the channel filtering, highlighting a stronger effect of the RIN on the upper eyes. Thus, when the system is limited by RIN rather than shot noise, the overall SNR is actually dominated by the lowest among the three SNRs of each of the three 2-PAM eyes, and the total BER can be approximated by the average among the three corresponding BERs. The RIN impact on the eye diagram is expected to become stronger for high ER values, thus in Section \ref{validity}A we will perform a detailed analysis of the SNR-to-BER conversion in extreme RIN conditions. In the remainder of the manuscript, unless otherwise stated, the global system BER will be calculated as the average of the per-eye BER. Fig. \ref{fig:model_snr}b shows the BER obtained in the same conditions as in Fig. \ref{fig:model_snr}a, compared to the BER computed in the simulations through error counting. The following equation relating SNR at the equalizer output and BER for a generic M-PAM modulation has been used to convert the SNR values returned by the analytical model for each of the three eyes:
\begin{equation}\label{eq:ber}
BER \simeq \frac{M-1}{M \cdot log_2(M)} \mbox{erfc} \left(\sqrt{\frac{3\cdot SNR}{2(M^2-1)}}\right)
\end{equation}
where $M$ is the constellation size. Through equation \ref{eq:ber} our model can be employed to estimate performance of any M-PAM modulation scheme, with any cardinality. Fig. \ref{fig:model_snr}b shows an excellent agreement of the analytical model with the time domain simulations, also in terms of BER, at least in the BER range of interest down to $10^{-4}$. The BER counting technique included in the simulator operates on about 500000 bits and can, therefore, become inaccurate for BER below $10^{-5}$. The inset in Fig. \ref{fig:model_snr}b show a perfect BER agreement also for FFE-based 25 GBaud 2-PAM and 8-PAM modulations with supergaussian filter order 1, confirming that the model can be applied to any M-PAM format. Similar results are obtained for DFE equalization as well (not shown due to space limitations). A pair of SNR values for the two equalization schemes can be obtain in about 0.045 seconds through the analytical model on a commercial standard PC. Compared to the 18.5 seconds needed to run the time domain simulation, the analytical tool thus yields over 410 times reduction in computational time consumption.
%
\section{Validity of the proposed model}
\label{validity}

In this section we check the validity of the model in different scenarios extending the analysis to somewhat "extreme" RIN conditions and including the effect of chromatic dispersion and APD-induced shot noise. Hereafter, the channel is described by a supergaussian filter of order 1.

\subsection{Effect of RIN}

\begin{figure}[h]
\centering
\includegraphics[width=1\linewidth]{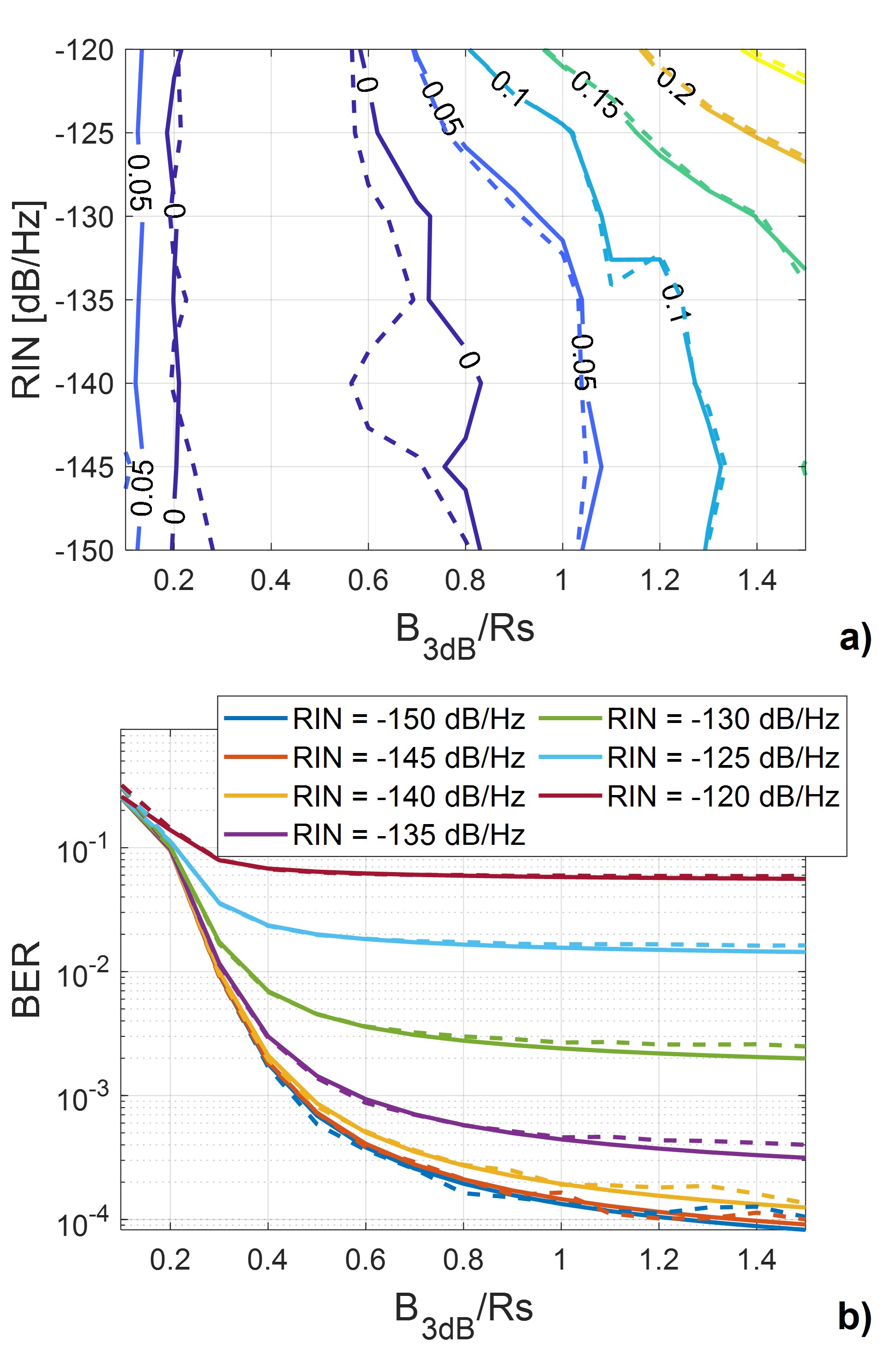}
\caption{a) $\Delta SNR$ difference in dB between the analytical and simulated results as a function of the RIN coefficient and of the ratio between the filter 3 dB bandwidth and the symbol rate for both FFE (solid) and DFE (dashed) equalization. b) BER evolution obtained through analytical model (solid) and time domain simulations (dashed) as a function of the ratio between the filter 3 dB bandwidth and the symbol rate for several values of the RIN coefficient and FFE equalization.}
\label{fig:rin}
\end{figure}

As mentioned in the previous Section \ref{model}, our model relies on the approximation that RIN and shot noise scale with the average useful signal power (squared for RIN). However, a more accurate model should assume that the variance is proportional to the instantaneous power. This approximation might become inaccurate when the noise levels are sufficiently high to break the validity of the assumption, especially for the RIN, as it is proportional to the square of the instantaneous signal power. Thus, we have performed a sort of stress test of the proposed model varying the RIN coefficient in the range from -150 $dB/Hz$ up to -120 $dB/Hz$ with ER = 12 dB and $P_{TX}$ = 0 dBm (OMA = 2.46 dBm). Moreover, to check the applicability of our model to any data rate we have increased the symbol rate of the simulated IMDD system to 56 GBaud, a typical value for modern 400G Ethernet solutions. The results of the comparison with the time domain simulator in terms of $\Delta SNR = SNR_{model}-SNR_{sim}$ in dB are depicted in Fig. \ref{fig:rin}a. The contour plot shows an estimation error very close to 0 dB across the whole investigated space of parameters for both FFE (solid) and DFE (dashed) scheme, with values above 0.05 dB only for very strongly bandlimited conditions ($B_{3dB}/R_s = 0.1$) or for filter bandwidth larger than the baud rate ($B_{3dB}/R_s > 1$). In this latter case a higher estimation error is to be expected: when there is no bandwidth limitation the need for channel equalization is reduced and so is the time averaging effect associated with the equalizer. The average power assumption our model is based upon is, thus, broken and the estimation becomes less accurate. Nevertheless, maximum errors of the order of 0.2 dB are still acceptable. Similarly, the evolution of the corresponding BER (obtained as the average of the three per-eye BER) in Fig. \ref{fig:rin}b shows a remarkable agreement between the two methods with FFE equalization, with negligible mismatch only for BER below $10^{-5}$, due to the limited simulated number of transmitted bits, and for $B_{3dB}/R_s$ ratios larger than 1, where the SNR estimation error slightly increases. Similar results are obtained with DFE equalization, not shown here due to space limitations.

\subsection{Chromatic Dispersion Approximation}

\begin{figure}[h]
\centering
\includegraphics[width=1\linewidth]{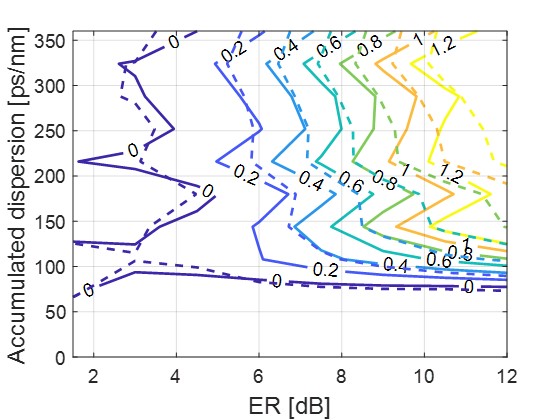}
\caption{$\Delta SNR$ difference in dB between the analytical and simulated results as a function of the accumulated dispersion and extinction ratio for both FFE (solid) and DFE (dashed) equalization.}
\label{fig:DSNR_CD}
\end{figure}

It is well known that in terms of field propagation the linear transfer function of a single mode fiber (SMF) in presence of chromatic dispersion only, can be described by the following expression:
\begin{equation}\label{eq:Hsmf}
H_{SMF}(f) = e^{j \pi c D L \cdot \frac{(f-f_c)^2}{f_c^2}}
\end{equation}
where $L$ is the fiber length, $f_c$ is the central frequency of the optical signal and $D$ is the CD coefficient. Unfortunately, an IMDD system is intrinsically nonlinear and thus the end-to-end electrical transfer function cannot be directly related to Eq. \ref{eq:Hsmf}. Nevertheless, in \cite{small_signal} a small signal analysis is carried out that allows to approximate Eq. \ref{eq:Hsmf} under the assumptions that the modulation amplitude is small compared to average signal power (i.e. that the M-PAM outer extinction ratio is small) and that the amplitude modulation is chirpless. The approximated small signal electrical to electrical transfer function of the system under chromatic dispersion effect only can be expressed as:
\begin{equation}\label{eq:Hcd}
H_{CD}(f) = cos\left[\pi cDL \left(\frac{f}{f_c}\right)^2\right]
\end{equation}
Including $H_{CD}(f)$ in our system, the $H_{ch}(f)$ function in Fig. \ref{fig:setup} and in Eq. \ref{eq:SNR_f_2} becomes now:
%

\begin{equation}\label{eq:Hchcd}
H_{ch}(f) = H'_{ch}(f)H_{CD}(f)
\end{equation}

where $H'_{ch}(f)$ is the channel frequency response without CD (i.e. the super-gaussian filter used in all our previous analyses). Note that the RIN contribution, as well as the useful signal, is now filtered also by the CD equivalent transfer function. Fig. \ref{fig:DSNR_CD} shows the comparison between the analytical model and the time domain simulations in terms of $\Delta SNR$ as a function of the accumulated dispersion ($D\cdot L$) and extinction ratio, when $RIN=-140$ dB/Hz, $P_{TX}=0$ dBm, $B_{3dB}/R_s=0.4$ and $\lambda_c = c/f_c = 1310$ nm.

The model estimation accuracy is excellent for accumulated dispersion up to 90 ps/nm and for ER up to 5 dB regardless of the considered equalization scheme. For a given value of the accumulated dispersion above 100 ps/nm the estimation error increases with ER, as the small signal approximation is gradually broken. The small signal approximation used in the derivation of Eq. \ref{eq:Hcd}, in fact, is based on the assumption that the amplitude modulation is small compared to the signal intensity, but increasing the ER results in a larger modulation amplitude. On the other hand, for a given ER an increase of the accumulated dispersion above 150 ps/nm does not alter the estimation accuracy substantially, indicating that the only limitation to our analytical model is the aforementioned small signal approximation.

\subsection{Avalanche Photodetection}

\begin{figure}[h]
\centering
\includegraphics[width=1\linewidth]{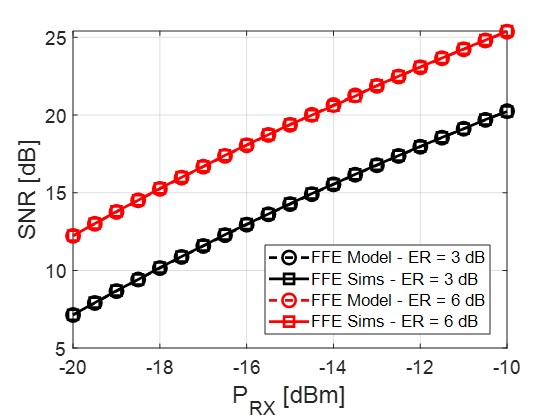}
\caption{SNR obtained through time domain simulations (solid, squares) and through the proposed analytical model (dashed, circles) as a function of the received optical power using 56 GBaud 4-PAM with FFE equalization. Transmitted power is 0 dBm and ER is 3 dB (black) or 6 dB (red).}
\label{fig:APD}
\end{figure}

In some applications where a higher receiver sensitivity is needed \cite{apd}, APDs are used to generate an electrical gain but also introducing an additional contribution to the total shot noise produced in the photodetection process. In fact, the PSD of the avalanche shot noise is given by:
\begin{equation}\label{eq:Pshot}
PSD_{shot} = \frac{G^2FqP_{RX}}{R}
\end{equation}
where $G$ is the gain of the photodiode, $F$ is its excess noise factor, $q$ is the electron charge and $R$ is the photodiode responsivity. In the following analysis we have considered APD parameters assuming 50G-class devices \cite{from25to50G,overview} with $G=10$ dB, $F=4.3$ dB, $R=0.7$ A/W and $B_{3dB}/R_s=0.7$. Fig. \ref{fig:APD} shows the SNR obtained with the two methods as a function of the received optical power for a 56 GBaud 4-PAM modulation when the transmitted power is 0 dBm and ER is 3 dB or 6 dB. Only results using FFE are shown as when the ratio $B_{3dB}/R_s$ is equal to 0.7 no significant SNR difference can be observed when DFE is used.

Fig. \ref{fig:APD} highlights again a nearly perfect match between the SNR estimated by the proposed analytical model and that computed through time domain simulations, regardless of the ER value. The estimation error is about 0.02 dB for each simulated $P_{RX}$ value.

\subsection{Equalizer Memory}

\begin{figure}[h]
\centering
\includegraphics[width=1\linewidth]{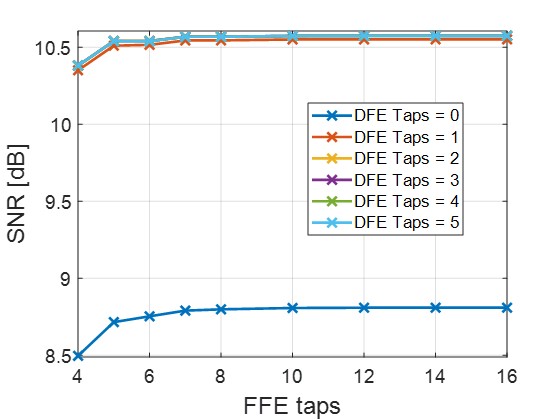}
\caption{SNR obtained through time domain simulations as a function of the number of FFE taps for several numbers of DFE taps. $P_{TX}=0$ dBm and $B_{3dB}/R_s=0.25$.}
\label{fig:taps}
\end{figure}

As the proposed analytical model is based on a derivation \cite{Fischer} that assumes infinitely long equalizer memory, in the previous sections we presented results obtained in comparison to time domain simulations performed using a high number of taps (see Table \ref{table_param}). Since equalizer convergence is usually sought with the shortest possible memory to reduce complexity and cost of the electronic implementation of the equalization scheme, here we investigate on the prediction accuracy of our model in comparison to more realistic versions of the equalization algorithm, varying the number of taps for both the FFE and DFE stage in the simulator. Fig. \ref{fig:taps} shows the SNR at the output of the equalizer as a function of the number of taps for a simulated 56 GBaud 4-PAM transmission with $P_{TX}=0$ dBm and a strong bandwidth limitation with $B_{3dB}/R_s=0.25$. When pure FFE equalization is used a number of taps as low as 10 is sufficient for perfect equalizer convergence and, even using 4 taps the SNR estimation error is only about 0.3 dB. Moreover, when DFE is used 10 taps for the FFE stage of the DFE equalizer are required to achieve optimum performance, whereas only 2 taps are enough for the DFE stage.
\section{Application Example}
\label{pon}

In this section we apply our analytical model to the analysis of the specific use case of a passive optical network. Standardization efforts for next generation PONs are being focused to the next big step in transmission speed \cite{overview}. Although it is yet not clear whether the next target bit rate will be 100G or 200G, and if the traditional IMDD-based PON architecture will transition to more complex yet advantageous coherent solutions \cite{scaling_laws}, in the mid term there will for sure be the need to upgrade current PON capacity to, at least, 100 Gbps/$\lambda$. In the PON environment this is usually intended as a gross target bit rate, thus including the overhead associated to the employed forward error correction (FEC) algorithm. The hard-decision FEC (HD-FEC) defined in the 50G PON standard sets the pre-FEC BER threshold to $10^{-2}$, but even soft-decision FECs (SD-FECs) might be selected in future PON standards raising the BER threshold to $1.9 \cdot 10^{-2}$ \cite{fec}. Moreover, current standards require at least 29 dB optical power budget (OPB) for N1 class PON. This physical layer requirement is hardly going to be modified and will be the main challenge in future versions of PON due to the reduced sensitivity of IMDD systems at such high bit rates.

\begin{figure}[h]
\centering
\includegraphics[width=1\linewidth]{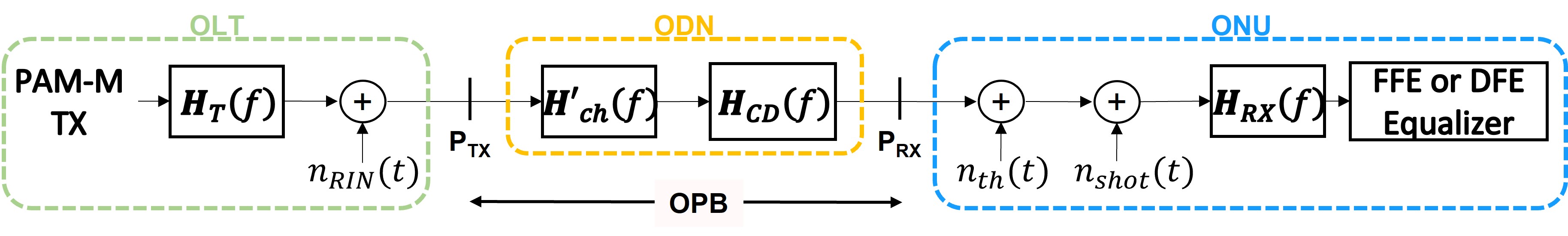}
\caption{Simplified scheme of the PON under investigation.}
\label{fig:pon_setup}
\end{figure}
Fig. \ref{fig:pon_setup} shows a simplified setup of a PON with APD-based detection, relating the optical line terminal (OLT), the optical distribution network (ODN) and the optical network unit (ONU) to the building blocks of our analytical model. As in the previous section we have considered APD parameters assuming 50G-class devices \cite{from25to50G,overview} with $G=10$ dB, $F=4.3$ dB and $R=0.7$ A/W. The OPB is defined across the ODN as the ratio between the transmitted optical power $P_{TX}$ and the received optical power $P_{RX}$ at a given target BER level. In our analysis we set the transmitted power to 11 dBm, a typical value often used in PONs to avoid the inset of nonlinear Kerr effects. Moreover, we assume the PON is operated in the O-band where the chromatic dispersion effect is small. Thus, in our approximation of the CD transfer function of the fiber we fixed the dispersion parameter $D$ to 3.85 $ps/(nm\cdot km)$, corresponding to the upper O-band wavelength of the grid plan of recently standardized 50G-PON \cite{overview}.
\begin{figure}[h]
\centering
\includegraphics[width=1\linewidth]{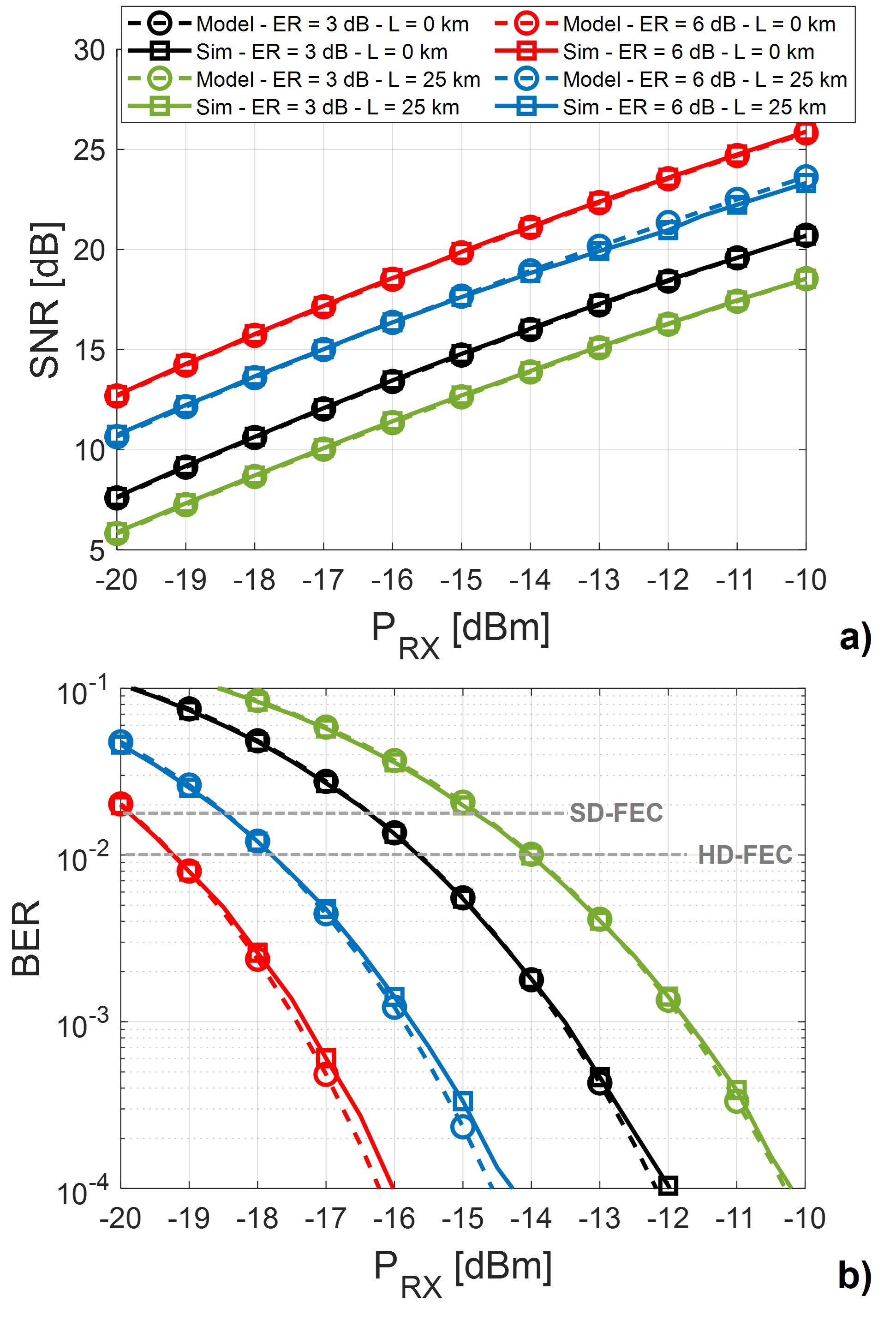}
\caption{a) SNR and b) corresponding BER obtained through time domain simulations (solid, squares) and through the proposed analytical model (dashed, circles) as a function of the received optical power using 50 GBaud 4-PAM (i.e. 100 Gbps) with FFE equalization in back-to-back (black, green) and with 25 km SMF in O-band (red, blue). Transmitted power is 11 dBm,  ER is 3 dB (black, red) or 6 dB (green, blue). Legend in a) applies to b) as well.}
\label{fig:pon_snr}
\end{figure}
Fig. \ref{fig:pon_snr} shows the comparison between model and simulations in terms of SNR and corresponding BER as a function of the received optical power for a 50 GBaud PON with two values of the ER, in back-to-back and with 25 km of single mode fiber. The ratio $B_{3dB}/R_s=0.7$, thus we only show results for FFE equalization. Similar results are obtained with DFE. Fig. \ref{fig:pon_snr}a highlights again a perfect agreement of the model with the simulations. A very small estimation error of about 0.3 dB on the SNR can be observed for received optical power above -14 dBm in the configuration with 25 km SMF and ER=6 dB, where the small signal approximation in Eq. \ref{eq:Hcd} starts to become inaccurate. In Fig. \ref{fig:pon_snr}b an excellent match on the sensitivity curves can also be observed. A slight mismatch of about 0.25 dB on the received optical power occurs at BER=$10^{-4}$ only for the two configurations with ER=6 dB. However, at the BER levels of interest around $10^{-2}$ the model is very accurate and shows two requirements for a PON with 25 km SMF and 100 Gbps 4-PAM modulation. To achieve the required 29 dB OPB the ER (typically the ER associated to a directly modulated laser or an electroabsorption modulated laser in the PON scenario) should be at least 6 dB and the use of a SD-FEC is necessary.
\section{Conclusions}
\label{conclusion}

We have presented an analytical model for fast performance estimation of IMDD-based communications systems affected by RIN, shot noise, thermal noise, chromatic dispersion and bandwidth limitations. The tool requires knowledge of the linear transfer function of the channel and can provide accurate SNR and BER prediction with over 410 times reduction in computation time compared to full time domain simulations.

We have evaluated the realm of applicability of our model focusing on 4-PAM modulation with FFE and DFE equalization for 50G and 100G target bit rates, showing SNR estimation errors below 0.1 dB for a wide range of RIN, shot noise and ER levels, regardless of the bandwidth limitation imposed by the channel. We also equipped our model with the possibility to take into account the effect of chromatic dispersion showing very good agreement with simulations for accumulated dispersion up to 360 $ps/(nm\cdot km)$ and ER up to 6 dB.

Lastly, we have applied and validated our analytical approach in the specific scenario of PAM-4 100 Gbps/$\lambda$ PON based on APD photodetection. Our findings reveal the need for at least 6 dB ER and advanced SD-FEC algorithms for the system to be able to meet the 29 dB OPB requirement imposed by current N1 PON class standard.

We believe that the proposed model is applicable to a wide range of current and future IMDD ultra high speed transmission and that it can be very useful, for instance, when a standardization body has to define a new transmission standard, and thus have to theoretically assess the expected performance in presence of different types of bandwidth limitations and noise effects. The model can be further extended to consider (but we cannot show the results here due to space limitations):

\begin{enumerate}
    \item unequally spaced M-PAM levels, for instance to study inaccuracies in M-PAM generation
    \item the small signal approximation we used in Eq. \ref{eq:Hcd} assumes a chirpless modulation and a linear fiber. There are models (again under the "small modulation" assumption) that enable the computation of an analytical transfer function also for transmitter chirp and nonlinear self-phase modulation along the fiber.
\end{enumerate}

The only impairment that our model cannot directly handle is the nonlinear distortion in the transmitted eye diagram associated, for instance, to the time-skew that can be present in directly modulated lasers \cite{reza, minelli} or in other nonlinear devices \cite{ramon}.

\section*{Acknowledgment}
This work was carried out under a research contract with Cisco Photonics and in the PhotoNext initiative at Politecnico di Torino. \url{www.photonext.polito.it}. Pablo Torres-Ferrera acknowledges the support of \emph{Secretaría de Educación, Ciencia, Tecnología e Innovación} (SECTEI), Mexico City. For the purpose of open access, the author has applied a Creative Commons Attribution (CC BY) license  to any Author Accepted Manuscript version arising. Data underlying the results presented here are available at https://doi.org/10.XXX/CAM.XXXX

\end{document}